\documentclass[graybox]{svmult}

\usepackage{mathptmx}       
\usepackage{helvet}         
\usepackage{courier}        
\usepackage{type1cm}        
\usepackage{makeidx}         
\usepackage{graphicx}        
\usepackage{multicol}        
\usepackage[bottom]{footmisc}

\usepackage{amsmath}
\usepackage{amsfonts}
\usepackage{xcolor}

\usepackage[utf8]{inputenc}

\usepackage{hyperref}
\hypersetup{colorlinks=true,linkcolor=blue,urlcolor=blue}

\def\argmin{\mathop{\text{argmin}}}
\def\mig{\frac{1}{2}}

\usepackage{ulem}
\definecolor{blue}{rgb}{0,0,1}

\long\def\pep#1#2{{\color{blue}\ifx#1\relax\else\sout{#1}\ \fi#2}}

\makeindex             

\begin{document}
	
	\title*{High Order Extrapolation Techniques for WENO Finite-Difference Schemes Applied to NACA Airfoil Profiles}
	\titlerunning{High Order Extrapolation for NACA Airfoil Profiles.}
	\author{Antonio Baeza, Pep Mulet and David Zor\'io}
	
	\institute{Antonio Baeza \at Departament de Matem\`atiques,  Universitat de Val\`encia, C/Dr. Moliner, 50, 46100, Burjassot, Val\`encia (Spain), \email{antonio.baeza@uv.es}
		\and Pep Mulet \at Departament de Matem\`atiques,  Universitat de Val\`encia, C/Dr. Moliner, 50, 46100, Burjassot, Val\`encia (Spain),  \email{mulet@uv.es} \and David Zor\'io \at Departament de Matem\`atiques,  Universitat de Val\`encia, C/Dr. Moliner, 50, 46100, Burjassot, Val\`encia (Spain),  \email{david.zorio@uv.es}}
	%
	
	%
	%
	
	\leavevmode\thispagestyle{empty}
	
	\noindent This version of the article has been accepted for publication, after a peer-review process, and is subject to Springer Nature’s AM terms of use, but is not the Version of Record and does not reflect post-acceptance improvements, or any corrections. The Version of Record is available online at: \url{https://doi.org/10.1007/978-3-319-63082-3_6}
	
	\newpage
	
	\maketitle
	
	\abstract{Finite-difference WENO schemes are capable of approximating accurately and efficiently weak solutions of hyperbolic
		conservation laws. In this context high order numerical boundary conditions
		have been proven to increase significantly the
		resolution of the 
		numerical solutions. In this paper a finite-difference WENO scheme
		is combined with a high order boundary extrapolation technique at
		ghost cells to solve problems involving NACA airfoil
		profiles. The results obtained are comparable with those obtained
		through other techniques involving unstructured meshes.}
	
	\section{Introduction}
	
	Cauchy problems for hyperbolic conservation laws usually involve
	weak solutions, which yields difficulties to tackle them numerically
	and has motivated the development of high resolution shock capturing
	schemes (onwards HRSC). In this paper we will focus on HRSC
	finite-difference schemes, that use a 
	discretization with a Cartesian mesh, to numerically solve problems
	involving NACA airfoil profiles, which are of a great interest in the
	field of aeronautical engineering.
	
	Due to the nature of a Cartesian mesh, if enough care is not taken at
	the boundary neighborhood the accuracy order will be lost. Therefore,
	a strategy to extrapolate the information from the computational
	domain plus the boundary conditions (if any) at the ghost cells should
	be developed, taking into account as well that singularities may be
	eventually positioned near the extrapolation zone in order to avoid an
	oscillatory behaviour of the scheme.
	
	Some authors have tackled this issue with different techniques, such
	as second order Lagrange interpolation with limiters in \cite{Sjogreen}
	and an inverse Lax-Wendroff procedure for
	inlet boundaries and high order least squares interpolation for outlet
	boundaries, equipped with scale and dimension dependent weights in
	order to account for discontinuities, in \cite{TanShu, TanWangShu}.
	
	The contents of the paper are organized as follows: Section
	\ref{sec:nsmp} stands for a brief description of some details about
	the numerical
	scheme and the procedure to automatically mesh complex domains through
	a Cartesian grid, as well as the criterion to choose the extrapolation
	stencil for the boundary conditions; in Section \ref{sec:extr} an
	extrapolation procedure with scale and dimension independent weights
	is described; in Section \ref{sec:ne} a
	numerical experiment involving a NACA0012 profile is shown as a matter
	of illustration; finally, in Section \ref{sec:scn} some conclusions
	are drawn.
	
	\section{Numerical Scheme and Meshing Procedure}\label{sec:nsmp}
	
	The equations that are considered here are $d$-dimensional $m\times m$ systems of
	hyperbolic conservation laws.
	\begin{equation}\label{eq:0}
		u_t+\nabla\cdot f(u)=0,\quad
		u:\Omega\times\mathbb{R}^+\subseteq\mathbb{R}^d\times\mathbb{R}^+\rightarrow\mathbb{R}^m,\quad
		f:\mathbb{R}^m\rightarrow\mathbb{R}^m,
	\end{equation}
	with prescribed initial condition $u(x;0)=u_0(x)$, $x\in\Omega$, and
	suitable boundary conditions.
	\subsection{Numerical Scheme}
	For the sake of simplicity, let us describe the procedure for $d=1$
	and $\Omega=(0,1)$. To approximate the solution of \eqref{eq:0}, these
	equations are then discretized in space by means of a Shu-Osher
	finite-difference approach \cite{ShuOsher1989}, by taking
	$x_j=(j+\mig)h$, $0\leq j<N$, $h=\frac{1}{N}$, therefore the spatial
	scheme can be written in local terms as follows:
	\begin{equation}\label{so}
		u_t=-\frac{\hat{f}_{j+\frac{1}{2}}(t)-\hat{f}_{j-\frac{1}{2}}(t)}{h}+\mathcal{O}(h^r),
	\end{equation}
	where
	$$\hat{f}_{j+\frac{1}{2}}(t)=\hat{f}(u(x_{j-k+1},
	t),\ldots,u(x_{j+k},t))$$
	are obtained through $r$-th order accurate WENO spatial
	reconstructions \cite{JiangShu96} 
	with the Donat-Marquina flux-splitting technique \cite{DonatMarquina96}.
	
	The ODE (\ref{so}) is then discretized (in time) through a Runge-Kutta
	procedure \cite{ShuOsher89}, yielding a fully discretized numerical
	scheme. This combination of techniques was proposed by Marquina and
	Mulet in \cite{MarquinaMulet03}.
	
	\subsection{Meshing Procedure}
	
	We summarize the meshing procedure for complex (non-rectandular) 
	domains. We focus on the two-dimensional case for simplicity.
	
	The computational domain is given by the intersections of the horizontal and vertical mesh lines
	that belong to the physical domain:
	$$\mathcal{D}:=\left((\overline{x}+h_x\mathbb{Z})\times
	(\overline{y}+h_y\mathbb{Z})\right)\cap\Omega,$$
	where $h_x,h_y>0$ are the grid sizes for each respective dimension and
	$\overline{x},\overline{y}$ are fixed real values.
	
	In order to advance in time using WENO schemes of order $2k-1$, 
	$k$ additional cells are needed at both sides of
	each horizontal and vertical mesh line. If these additional
	cells fall outside the domain they are usually named \textit{ghost
		cells}.
	
	A normal direction associated to each ghost cell $P$ and the domain
	boundary, $\partial\Omega$, is then computed by finding the
	point $P_0$ such that
	$$\|P-P_0\|_2=\min\{\|P-B\|_2:\quad B\in\partial\Omega\}$$
	and then considering the corresponding (interior) normal vector,
	$\overrightarrow{n}(P)=\overrightarrow{PP_0}$.
	
	To extrapolate the information at each ghost cell, we first
	interpolate/extrapolate information from the computational domain to
	points on the normal line, and then extrapolate the obtained
	data together with the boundary condition, if any, in order to fill
	the missing value of the corresponding ghost cell.
	
	\begin{figure}[htb]
		\centering
		\begin{tabular}{cc}
			\includegraphics[width=0.47\textwidth,height=4cm]{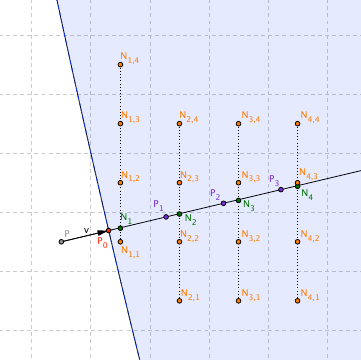} &
			\includegraphics[width=0.47\textwidth,height=4cm]{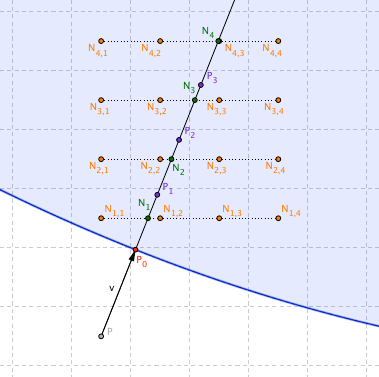}\\
			(a) & (b)
		\end{tabular}
		\caption{Examples of choice of stencil ($R=3$). We use the stencil
			$\mathcal{S}(P)=\{N_1, N_2, N_3, N_4\}$ in case of outflow
			boundary and conditions and the stencil $\mathcal{S}(P)=\{P_0,
			P_1, P_2, P_3\}$ in case of Dirichlet boundary conditions.}
		\label{fig:CxyW}
	\end{figure}
	
	Figure \ref{fig:CxyW} illustrates graphically the selection procedure
	of nodes from the computational domain, which can be described through
	the following summarized algorithm:
	
	\begin{enumerate}
		\item Extrapolate/interpolate the selected data from the
		computational domain, $N_{i,j}$, $1\leq i,j\leq R+1$, to the
		corresponding nodes $N_i$ on the normal line, $1\leq i\leq R+1$. The
		axis direction is chosen in terms of the angle of the normal vector.
		\item If there are prescribed boundary conditions, the
		information contained in the normal line involving the
		recently extrapolated/interpolated values at the points $N_i$,
		$1\leq i\leq R+1$,
		is used to interpolate/extrapolate them at new nodes $P_i$, $1\leq
		i\leq R$, in
		a way such that the stencil including the boundary node $P_0$
		(corresponding to the intersection of the normal line with
		$\partial\Omega$) are equally spaced.
		\item There are two possibilities:
		\begin{enumerate}
			\item If boundary conditions are prescribed, we extrapolate the
			information contained in the stencil $P_i$, $0\leq i\leq R$, at
			the ghost cell $P$.
			\item Otherwise, we extrapolate directly the information
			contained in the points obtained at the first step, $N_i$, $1\leq
			i\leq R+1$, to the ghost cell $P$.
		\end{enumerate}
	\end{enumerate}
	
	In case of reflecting boundary conditions for the two dimensional
	Euler equations, with velocity component $\overrightarrow{v}$, then
	one defines
	$$\overrightarrow{n}=\frac{\overrightarrow{P_0P}}{\big\Vert
		\overrightarrow{P_0P}\big\Vert},\quad
	\overrightarrow{t}=\overrightarrow{n}^{\perp},
	$$
	and obtains normal
	and tangential components of $\overrightarrow{n}$
	at each point $P_i$ of the mentioned segment by:
	\begin{equation*}
		v^t(P_i)=\overrightarrow{v}(P_i)\cdot \overrightarrow{t}, \quad
		v^n(P_i)=\overrightarrow{v}(P_i)\cdot \overrightarrow{n}.
	\end{equation*}
	The extrapolation procedure is applied to $v^t(P_i)$ to approximate
	$v^t(P)$  and to $v^n(P_i)$  and $v^n(P_0)=0$ to approximate
	$v^n(P)$. Once $v^t(P), v^{n}(P)$ are approximated, the
	approximation to $\overrightarrow{v}(P)$ is set to
	\begin{equation*}
		\overrightarrow{v}(P)=v^{t}(P)\overrightarrow{t}+v^{n}(P)\overrightarrow{n}.
	\end{equation*}
	
	More detailed information about the meshing procedure can be
	found at \cite{BaezaMuletZorio2015,BaezaMuletZorio2016}.
	
	\section{Extrapolation Procedure}\label{sec:extr}
	
	Likewise the WENO scheme idea, we need to ensure that the
	extrapolation at the ghost cells is performed smoothly, namely, either there
	are no discontinuities in the extrapolation stencil or, in such case,
	a procedure to detect them and take only information from the correct
	side of the discontinuity is developed.
	
	Let us consider a stencil of equally spaced nodes $x_0 < \dots < x_R$,
	their corresponding nodal values $u_j=u(x_{j})$ and $x_*$ the
	extrapolation node. Then we define
	$j_0$ as the index corresponding to the closest node in the
	stencil to the extrapolation node, namely
	$j_0=\argmin_{j\in \{0, \dots, R\}}|x_j-x_*|.$
	We then define the corresponding smoothness indicator for each stencil
	$\{x_k,\ldots,x_{k+r_0}\}$, $0\leq k\leq R-r_0$:
	$$I_k=\frac{1}{r}\sum_{\ell=1}^{r_0}\int_{x_0}^{x_r}h^{2\ell-1}p_k^{(\ell)}(x)^2dx,\quad
	0\leq k\leq R-r_0,$$
	where $p_k$ is the polynomial of degree at most $k$ such that
	$p_k(x_{i+j})=u_{i+j}$ for $0\leq j\leq r_0,$ $0\leq i\leq R-r_0$.
	
	If $r$-th order accuracy is desired, $r\leq R$, by stability
	motivations (see \cite{BaezaMuletZorio2016}), we now define $p$ to be a
	polynomial of degree $r$ satisfying $p(x_i)=u_i$, $0\leq i\leq R$, by
	least squares.
	
	Finally, we account for discontinuities by using properly the
	smoothness indicators computed previously. We thus define the
	following global average weight
	$$\omega:=\frac{(R-r_0+1)^2}{\left(\sum_{k=0}^{R-r_0}I_k\right)
		\left(\sum_{k=0}^{R-r_0}\frac{1}{I_k}\right)}.$$
	This weight satisfies $0\leq\omega\leq1$ and
	$$\omega=\left\{\begin{array}{rl}
		1-\mathcal{O}(h^{2c_s}) & \textnormal{if the stencil is
		}C^{r_0}\textnormal{ with a }s\textnormal{-th order zero derivative}, \\
		\mathcal{O}(h^{2}) & \textnormal{if the stencil contains a
			discontinuity,}
	\end{array}\right.$$
	with $c_s:=\max\{r_0-s,0\}.$
	
	The final result of the extrapolation is thus defined by
	$u_*=\omega p(x_*)+(1-\omega) u_{i_0}.$
	The above expression can be interpreted as a continuous transition
	between the full order extrapolation in case of smoothness in the
	stencil ($\omega\approx1$) and the reduction to a first order constant
	extrapolation, taking the closest node from the stencil with respect
	to the extrapolation point as the reference value, in presence of
	discontinuities on the stencil ($\omega\approx0$).
	\section{Numerical Experiment}\label{sec:ne}
	
	In this section we show a numerical experiment with the flow given by
	the 2D Euler equations (with adiabatic constant $\gamma=1.4$) around a NACA
	airfoil profile for $R=8$, $r=4$, $r_0=2$, WENO5 spatial discretization and RK3
	time scheme (fifth order spatial accuracy and third order time
	accuracy). Prior tests in \cite{BaezaMuletZorio2016} showed that this
	scheme is indeed fifth order accurate in space in smooth problems,
	even when complex geometries are considered.
	
\begin{figure}[htb]
		\centering
		\begin{tabular}{cc}
			\includegraphics[height=0.4\textheight]{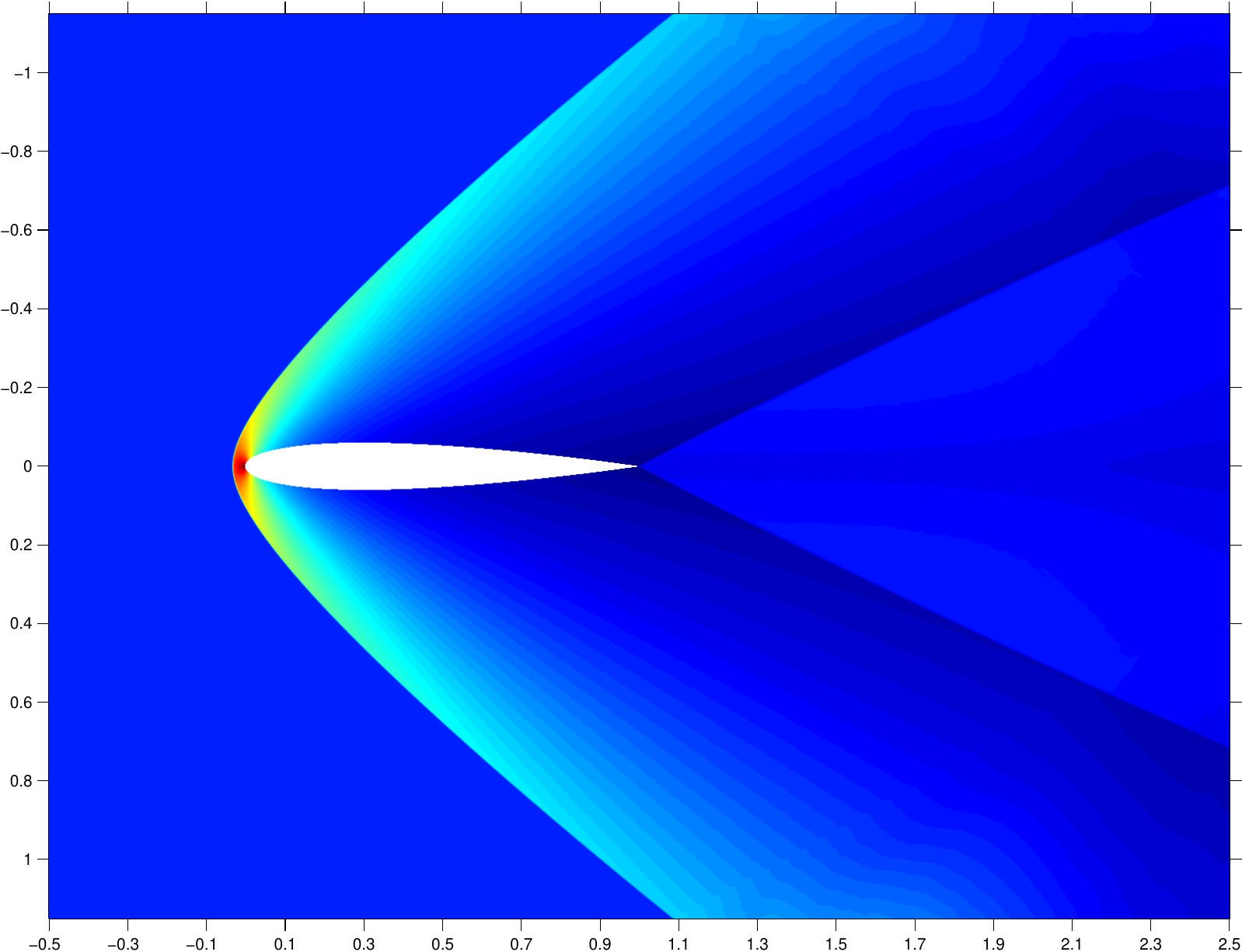} &
			\includegraphics[height=0.4\textheight]{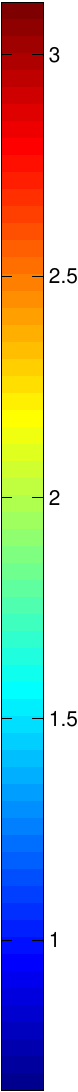} \\
			(a) Density field\\
			\includegraphics[height=0.4\textheight]{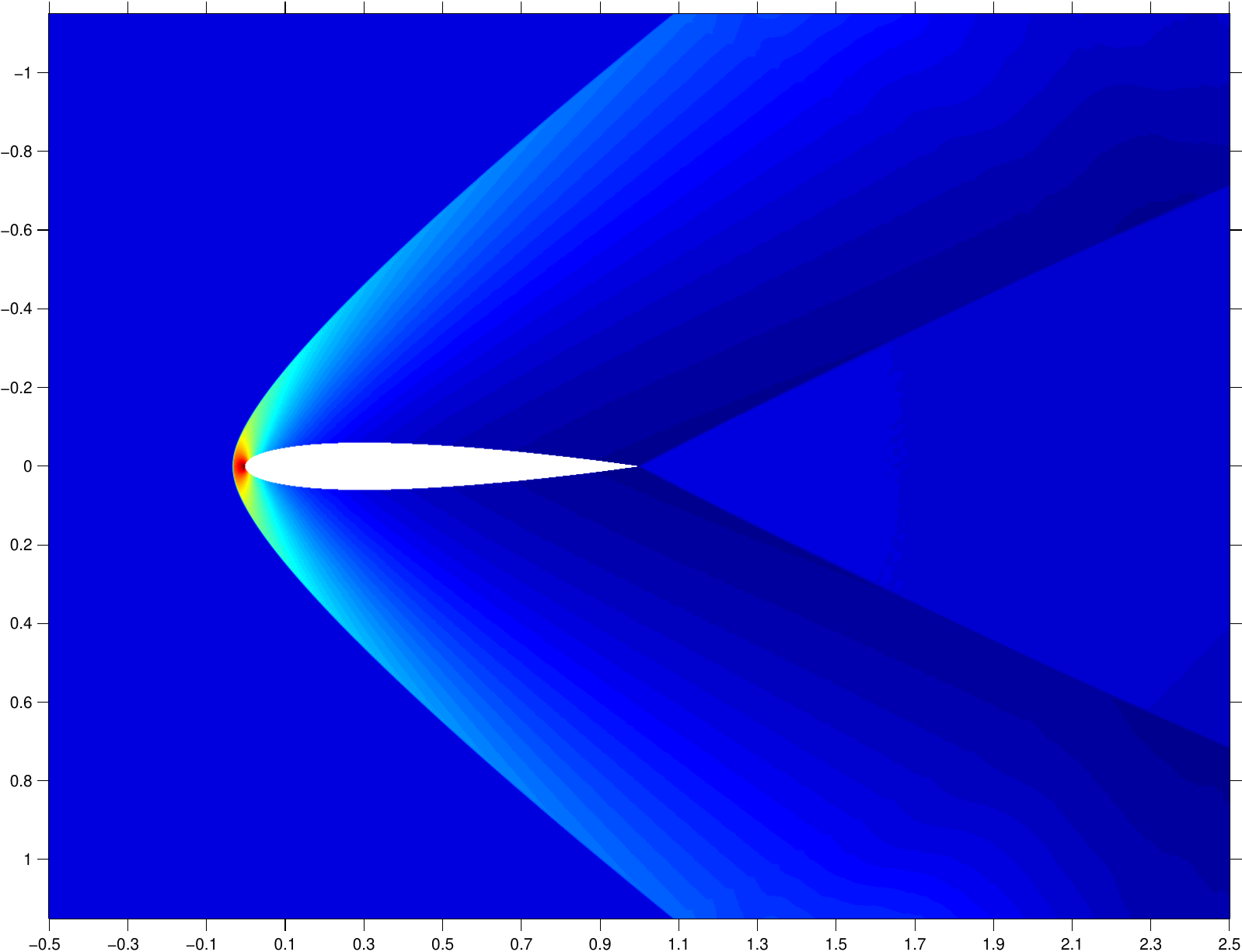} &
			\includegraphics[height=0.4\textheight]{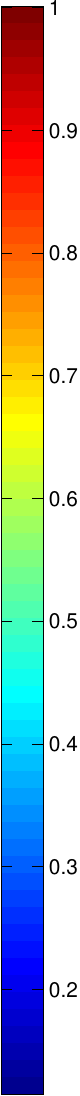} \\
			(b) Pressure field
		\end{tabular}
		\caption{Density and pressure fields for the supersonic NACA0012 interaction
			test.}
		\label{naca}
	\end{figure}

	We consider a Mach 2 flow interacting with a NACA0012 profile
	\cite{Feistauer} and we
	perform the simulation until $T=5$, whose results are shown in Figure
	\ref{naca} for the density and pressure fields.
	
	The results are consistent with those obtained in Feistauer et
	al. \cite{Feistauer} using finite volume techniques.
	\section{Conclusions}
	In this paper a technique to perform high order boundary
	extrapolations in finite difference methods has been
	successfully applied to supersonic inviscid flow around a slender NACA
	airfoil.
	
	The extrapolation attains the desired accuracy in smooth problems as
	shown in \cite{BaezaMuletZorio2016} and behaves properly in presence
	of discontinuities, which is useful for problems of practical interest
	like the one presented herein.
	\label{sec:scn}
	
	%
	%
	%
	\begin{acknowledgement}
		This research was partially supported by Spanish MINECO grant MTM2014-54388-P.
	\end{acknowledgement}

\end{document}